\documentclass[reqno,tbtags,a4paper,12pt]{amsart}

\usepackage[in]{fullpage}

\usepackage{amsmath,amssymb,amsthm,amsfonts,amscd}

\newtheorem*{theorem*}{Theorem}

\newtheorem*{cor*}{Corollary}

\newcommand{\Z}{\mathbb{Z}}
\newcommand{\R}{\mathbb{R}}
\newcommand{\M}{\mathcal{M}}
\newcommand{\s}{\mathcal{S}}
\renewcommand{\P}{\mathcal{P}}

\author{Alexander A. Gaifullin}
\thanks{The work was partially supported by the Russian Foundation
for Basic Research (grant no.~06-01-72551) and the Russian
Programme for the Support of Leading Scientific Schools (grant
no.~1824.2008.1).}

\title{Realisation of cycles by aspherical manifolds}
\date{}

\address{Moscow State University}

\begin{document}
\maketitle

In the late 1940s N.~Steenrod posed the following problem, which
is now familiar as the problem on realisation of cycles. For a
given homology class $z\in H_n(X;\Z)$, do there exist an oriented
manifold~$N^n$ and a mapping $f:N^n\to X$ such that $f_*[N^n]=z$?
A famous theorem of R.~Thom claims that each integral homology
class is realisable in sense of Steenrod with some multiplicity. A
classical problem is the problem of realisation of cycles by
images of spheres, that is, the problem of the description for the
image of the Hurewicz homomorphism. In this case not every
homology class can be realised with multiplicity. It is
interesting to find a class~$\M_n$ of smooth $n$-dimensional
manifolds sufficient for realisation with multiplicities of all
integral $n$-dimensional homology classes of every space~$X$. The
following theorem is the main result of this paper.

\begin{theorem*}
Suppose $M^n$ is the isospectral manifold of real symmetric
tridiagonal $(n+1)\times (n+1)$ matrices and $X$ is an arbitrary
arcwise connected topological space. Then for each homology class
$z\in H_n(X;\Z)$ there are a connected finite-fold covering
$p:\widehat{M}^n\to M^n$ and a mapping $f:\widehat{M}^n\to X$,
such that $f_*[\widehat{M}^n]=qz$ for some positive integer~$q$.
\end{theorem*}

C.~Tomei~\cite{Tom84} proved that~$M^n$ is an aspherical smooth
oriented manifold. (A manifold is called {\it aspherical} if it
has homotopy type $K(\pi,1)$.) The group $\pi_1(M^n)$ was computed
by M.~Davis~\cite{Dav87}. He proved that it is isomorphic to a
torsion-free subgroup of finite index of the Coxeter group
\begin{multline*}
W=\langle s_1,\ldots,s_n,r_1,\ldots,r_n\mid s_i^2=r_i^2=1,
s_is_j=s_js_i\text{ for
}|i-j|>1,\\s_is_{i+1}s_i=s_{i+1}s_is_{i+1},
r_ir_j=r_jr_i,s_ir_j=r_js_i\text{ for }i\ne j\rangle.
\end{multline*}
\begin{cor*}
Every integral homology class of every arcwise connected space can
be realised  with some multiplicity by an image of an oriented
aspherical smooth manifold with fundamental group isomorphic to a
torsion-free subgroup of finite index of the group~$W$.
\end{cor*}

Put $[n+1]=\{1,2,\ldots,n+1\}$. Let $\s$ be the set of all
nonempty subsets $\omega\subset[n+1]$, $\omega\ne[n+1]$. The {\it
permutahedron}~$\Pi^n$ is the convex hull of the points obtained
by all possible permutations of coordinates of the point
$(1,2,\ldots,n+1)\in\R^{n+1}$. The permutahedron is an
$n$-dimensional simple convex polytope. Its facets are in
one-to-one correspondence with subsets~$\omega\in\s$. Facets
$F_{\omega_1}$ and $F_{\omega_2}$ intersect each other if and only
if either $\omega_1\subset\omega_2$ or $\omega_2\subset\omega_1$.
Any face of the permutahedron~$\Pi^n$ has form
$F_{\omega_1}\cap\ldots\cap F_{\omega_k}$, where
$\varnothing\varsubsetneq\omega_1\varsubsetneq\ldots\varsubsetneq\omega_k\varsubsetneq[n+1]$;
the barycenter of this face will be denoted by
$b_{\omega_1,\ldots,\omega_k}(\Pi^n)$.

C.~Tomei~\cite{Tom84} constructed a decomposition of the
manifold~$M^n$ into $2^n$ permutahedra of the form
$M^n=(\Z_2^n\times\Pi^n)/\sim$, where the equivalence
relation~$\sim$ is generated by the identifications $(g,x)\sim
(e_{|\omega|}g,x)$ for $x\in F_{\omega}$. Here $e_1,\ldots,e_n$
are the generators of the group~$\Z_2^n$. Let $[g,x]$ be the
equivalence class of~$(g,x)$. Let us give an explicit construction
of a decomposition of the manifold $\widehat{M}^n$ into
permutahedra covering the decomposition of the manifold~$M^n$.
Every homology class of an arcwise connected space can be realised
by an image of a strongly connected oriented pseudo-manifold~$Z^n$
(for definition, see~\cite[\S 24]{SeTh38}). Hence we suffice to
prove Theorem only for the case $X=Z^n$ and $z=[Z^n]$. Taking the
first barycentric subdivision of the complex~$Z^n$ we may assume
that its vertices are regularly coloured in colours from  the
set~$[n+1]$. (``Regularly'' means that any two vertices connected
by an edge are coloured in distinct colours.) By~$\mu(\sigma)$ we
denote the set of colours of vertices of a simplex~$\sigma$. For
an $n$-dimensional simplex~$\sigma$ by $b_{\omega}(\sigma)$ we
denote the barycenter of the face $\tau\subset\sigma$ such that
$\mu(\tau)=\omega$.

By $U$ we denote the set of $n$-dimensional simplices of the
complex~$Z^n$. Since the vertices of~$Z^n$ admit a regular
colouring, we obtain that the set~$U$ can be decomposed into two
parts $U_+\sqcup U_-$ so that two simplices possessing a common
facet belong to distinct parts. For each $\omega\in\s$ we denote
by~$\P_{\omega}$ the set of involutions $\Lambda:U\to U$ such that
$\Lambda(U_{\pm})=U_{\mp}$ and
$\mu\left(\sigma\cap\Lambda(\sigma)\right)\supset\omega$ for every
simplex~$\sigma\in U$. The sets~$\P_{\omega}$ are nonempty. Define
a homomorphism $\eta:\Z_2^n\to\Z_2$ on the generators
by~$\eta(e_i)=-1$. Define a set~$V$ and
involutions~$\Phi_{\omega}:V\to V$ by
\begin{gather*}
V=\left(U_+\times\prod_{\omega\in\s}\P_{\omega}\times\eta^{-1}(1)\right)\cup
\left(U_-\times\prod_{\omega\in\s}\P_{\omega}\times\eta^{-1}(-1)\right)\subset
U\times\prod_{\omega\in\s}\P_{\omega}\times\Z_2^n;\\
\Phi_{\omega}\left(\sigma,\left(\Lambda_{\gamma}\right)_{\gamma\in\s},g\right)=
\left(\Lambda_{\omega}(\sigma),\left(\widetilde{\Lambda}_{\gamma}\right)_{\gamma\in\s},
e_{|\omega|}g\right),
\end{gather*}
where $\widetilde{\Lambda}_{\gamma}=\Lambda_{\omega}\circ
\Lambda_{\gamma}\circ\Lambda_{\omega}$ if $\gamma\subset\omega$
and $\widetilde{\Lambda}_{\gamma}=\Lambda_{\gamma}$ if
$\gamma\not\subset\omega$.

Put $\widehat{M}^n=V\times\Pi^n/\sim$, where the equivalence
relation~$\sim$ is generated by the identifications~$(v,x)\sim
(\Phi_{\omega}(v),x)$ for $x\in F_{\omega}$.

The projection $p:\widehat{M}^n\to M^n$ is given by
$p([(\sigma,(\Lambda_{\omega})_{\omega\in\s},g),x])=[g,x]$. The
mapping~$p$ is well defined and is a finite-fold covering.
Therefore $\widehat{M}^n$ is a smooth oriented manifold.

Let $K$ be the barycentric subdivision of the constructed
decomposition of the manifold~$\widehat{M}^n$ into permutahedra.
Define a mapping $f:\widehat{M}^n\to Z^n$ on the vertices of the
triangulation~$K$ by
$f([v,b_{\omega_1,\ldots,\omega_k}(\Pi^n)])=b_{\omega_1}(\sigma)$
and extend it linearly to every simplex of~$K$. The mapping~$f$ is
well defined and $f_*[\widehat{M}^n]=q[Z^n]$, where
$q=2^{n-1}\prod_{\omega\in\s}|\P_{\omega}|$.

The obtained manifold~$\widehat{M}^n$ is not necessarily
connected. The required manifold is an arbitrary connected
component of the manifold~$\widehat{M}^n$.

The manifold~$\widehat{M}^n$ implicitly appeared in the author's
paper~\cite{Gai07} as a special case of a general construction of
realisation of cycles. The decomposition of the
manifold~$\widehat{M}^n$ constructed in the present paper is
distinct from the decomposition constructed in~\cite{Gai07} and is
adapted for proving that $\widehat{M}^n$ covers~$M^n$.

The author is grateful to V.M.~Buchstaber for posing the problem
and permanent attention and to S.M.~Natanzon, A.V.~Penskoi,
A.B.~Sosinsky and O.V.~Schwarzman for useful discussions.

\end{document}